\documentclass[12pt]{article}
\usepackage{amssymb,amsmath,mathrsfs}
\usepackage{hyperref}
\usepackage{graphicx}
\usepackage{color}
\usepackage[OT2,OT1]{fontenc}
\usepackage{upquote}
\usepackage[numbers,sort&compress]{natbib}
\usepackage{alphabeta}

\begin{document}

\title{\Large\bf An iterative method for computing $\pi$ by argument reduction of the tangent function}

\author{
\normalsize\bf Sanjar M. Abrarov, Rehan Siddiqui, Rajinder Kumar Jagpal \\
\normalsize\bf and Brendan M. Quine}

\date{February 28, 2024}
\maketitle

\begin{abstract}
In this work, we develop a new iterative method for computing the digits of $\pi$ by argument reduction of the tangent function. This method combines a modified version of the iterative formula for $\pi$ with squared convergence that we proposed in a previous work and a leading arctangent term from the Machin-like formula. The computational test we performed shows that algorithmic implementation can provide more than $17$ digits of $\pi$ per increment. Mathematica codes, showing the convergence rate for computing the digits of $\pi$, are presented.
\vspace{0.2cm}
\\
\end{abstract}

\section{Introduction}

In 1706, English astronomer and mathematician John Machin discovered a formula for $\pi$, expressed in terms of two arctangents:
\begin{equation}\label{OMF} 
\frac{\pi}{4} = 4\arctan\left(\frac{1}{5}\right) - \arctan\left(\frac{1}{239}\right).
\end{equation}
By using this formula, he was the first to compute $100$ digits of $\pi$ \cite{Beckmann1971,Berggren2004,Borwein2008,Agarwal2013}. Nowadays,~identities of the kind
\begin{equation}\label{GFM} 
\frac{\pi}{4} = \sum\limits_{j = 1}^J\alpha_j\arctan\left(\frac{1}\beta_j\right),
\end{equation}
where the coefficients $\alpha_j$ and $\beta_j$ are either integers or rational numbers, are named after him as the Machin-like formulas for $\pi$.

The Maclaurin series expansion of the arctangent function is given by
\begin{equation}\label{MSE4AF} 
\arctan\left(x\right) = x - \frac{x^3}{3!} + \frac{x^5}{5!} - \frac{x^7}{7!} + \cdots = \sum\limits_{n=0}^\infty \frac{(-1)^n x^{2n+1}}{2n + 1},
\end{equation}
from which it follows that
\begin{equation}\label{ABON} 
\arctan(x) = x + O(x^3).
\end{equation}
Consequently, we can conclude that the convergence of the Machin-like Formula \eqref{GFM} for $\pi$ is better when the coefficients $\beta_j$ are larger by their absolute~values. 

In order to estimate the efficiency of a given Machin-like formula, Lehmer introduced the measure~\cite{Lehmer1938,Castellanos1988,Tweddle1991,Abeles1993,Wetherfield1996}:
\begin{equation}\label{LM} 
\mu = \sum\limits_{j = 1}^J \frac{1}{\log_{10}\left(|\beta_j|\right)}.
\end{equation}

A smaller value of $\mu$ indicates a higher efficiency of a given Machin-like formula. According to this formula, the~measure $\mu$ is smaller when the total number of the terms $J$ is lower and the coefficients $\beta_j$ are higher by their absolute values. A~more detailed description and the significance of the measure \eqref{LM} can be found in the literature~\cite{Castellanos1988,Wetherfield1996}.

Although Lehmer's paper~\cite{Lehmer1938} assumes that $\beta_j\in\mathbb{Q}$, it seems that the measure \eqref{LM} is applicable only when all values of $\beta_j$ are integers. Otherwise, even if a single value from the set $\left\{\beta_j\right\}$ represents a quotient, its presence causes complexities in the computation~\cite{Gasull2023}, which appear due to the rapidly increasing number of digits in the numerators in summation terms when a series expansion of the arctangent function, like Equations~\eqref{MSE4AF}, \eqref{EF} or \eqref{AFSE1}, is used in the computation. That is why it is very desirable to generate the Machin-like formulas for $\pi$ where $\beta_j\in\mathbb{Z}$.

The following Machin-like formula~\cite{Lehmer1938}
\[
\frac{\pi}{4} = 12\arctan\left(\frac{1}{18}\right) + 8\arctan\left(\frac{1}{57}\right) - 5\arctan\left(\frac{1}{239}\right),
\]
discovered by Gauss, has Lehmer's measure $\mu \approx 1.78661$. This signifies that this equation is more efficient in computing digits of $\pi$ since its measure is smaller than the original Machin-like Formula~\eqref{OMF} with measure $\mu \approx 1.85113$.

All coefficients in any Machin-like formula for $\pi$ satisfy the relation
\begin{equation}\label{PR} 
\prod\limits_{j = 1}^J (\beta_j + i)^{\alpha_j}\propto(1 + i)
\end{equation}
implying that the real and imaginary parts of this product must be equal to each other. This relation can be used for validation. For~example, there is a simple and elegant proof of the original Machin Formula~\eqref{OMF} for $\pi$ \cite{Guillera2009}
\[
(5 + i)^{4}(259 + i)^{-1} = 2(1 + i).
\]

Although several iterative formulas with quadratic, cubic, quartic, quintic and nonic convergences have been discovered~\cite{Borwein1984,Borwein1986a,Borwein1986b,Borwein1987,Bailey1997}, they require undesirable surd numbers that appear over and over again at each consecutive step of~the iteration.

Historically, algorithms based on the Chudnovsky and Machin-like formulas successively broke  records in computing the digits of $\pi$ \cite{Agarwal2013}. Currently, the~Chudnovsky formula, providing linear convergence with $14$ to $16$ digits of $\pi$ per increment of the summation terms~\cite{Beckmann1971,Berggren2004,Agarwal2013},  appears to be the most~efficient.

In 2002, Kanada broke a record by computing more than one trillion decimal digits of $\pi$ for the first time by using the following self-checking pair (also known as the St\"{o}rmer--Takano pair) of Machin-like formulas~\cite{Calcut2009,Agarwal2013}
\[
\begin{aligned}
\frac{\pi }{4} &= 44\arctan \left(\frac{1}{57}\right)+7\arctan\left(\frac{1}{239}\right) - 12\arctan\left(\frac{1}{682}\right) \\ 
&+24\arctan\left(\frac{1}{{12943}}\right)  
\end{aligned}
\]
and
\[
\begin{aligned}
\frac{\pi}{4} &= 12\arctan\left(\frac{1}{49}\right) + 32\arctan\left(\frac{1}{57}\right) - 5\arctan\left(\frac{1}{239}\right) \\ 
&+12\arctan\left(\frac{1}{110443}\right)  
\end{aligned}
\]
with Lehmer's measures of $1.58604$ and $1.7799$, respectively. Although~the current record, achieved by using the Chudnovsky formula, exceeds one hundred trillion digits of $\pi$ \cite{Gasull2023}, application of the Machin-like formulas may be promising to calculate a comparable number of digits due to availability of more advanced and powerful supercomputers then those used by Kanada more than $20$ years ago. Furthermore, more Machin-like formulas with smaller Lehmer's measures have been  discovered~\cite{Chien-Lih1997,Nimbran2010,Wetherfield2013}. For~example, the~following two equations
\begin{equation}\label{SChP1} 
\begin{aligned}
\frac{\pi}{4} &= 83\arctan \left(\frac{1}{107}\right) + 17\arctan\left(\frac{1}{1710}\right) - 22\arctan\left(\frac{1}{103697}\right) \\
&-24\arctan\left(\frac{1}{2513489}\right) - 44\arctan\left(\frac{1}{18280007883}\right) \\ 
&+12\arctan\left(\frac{1}{7939642926390344818}\right) \\ 
&+22\arctan\left(\frac{1}{3054211727257704725384731479018}\right)
\end{aligned}
\end{equation}
and
\begin{equation}\label{SChP2} 
\begin{aligned}
\frac{\pi}{4} &= 83\arctan\left(\frac{1}{107}\right) + 17\arctan\left(\frac{1}{1710}\right) - 22\arctan\left(\frac{1}{103097}\right) \\ 
&-12\arctan\left(\frac{1}{1256744}\right) - 22\arctan\left(\frac{1}{9140003941}\right) \\ 
&+12\arctan\left(\frac{1}{3158812219818}\right) \\ 
&+22\arctan\left(\frac{1}{167079344092131066905}\right),
\end{aligned}
\end{equation}
may be more efficiently used as a self-checking pair since their Lehmer's measures of $1.34085$ and $1.39524$, respectively, are considerably smaller. Therefore, application of  Machin-like formulas with small Lehmer's measure have colossal potential and can be competitive for computing the digits of the constant $\pi$.

Equation \eqref{SChP1} was reported by Wetherfield~\cite{Wetherfield2013}. Equation \eqref{SChP2} was obtained from another Wetherfield identity~\cite{Wetherfield2013}
\[
\begin{aligned}
\frac{\pi}{4} &= 83\arctan\left(\frac{1}{107}\right) + 17\arctan\left(\frac{1}{1710}\right) - 22\arctan\left(\frac{1}{103697}\right) \\ 
&-12\arctan\left(\frac{2}{2513489}\right) - 22\arctan\left(\frac{2}{18280007883}\right).  
\end{aligned}
\]
by using the identity~\cite{Abrarov2022}
\[
\arctan\left(\frac{1}{z}\right) = \arctan\left(\frac{1}{\left\lfloor z \right\rfloor}\right) + \arctan\left(\frac{\left\lfloor z \right\rfloor - z}{1 + z\left\lfloor z \right\rfloor}\right), \qquad z \notin\left[0,1\right).
\]

In our earlier publication~\cite{Abrarov2017}, we derived the following two-term Machin-like formula
\begin{equation}\label{TTMF1} 
\frac{\pi}{4} = 2^{k-1}\arctan\left(\frac{1}{\gamma}\right) + \arctan\left(\frac{1 - \sin\left(2^{k-1}\arctan\left(\frac{2\gamma}{\gamma^2 - 1}\right)\right)}{\cos\left(2^{k-1}\arctan\left(\frac{2\gamma}{\gamma^2 - 1}\right)\right)}\right),
\end{equation}
where the constant $\gamma$ may be conveniently chosen according to a relation~\cite{WolframCloud}
\begin{equation}\label{GA} 
\frac{2^{k - 1}}{\gamma} \approx \frac{\pi}{4}.
\end{equation}
The complete proof of the identity \eqref{TTMF1} is quite lengthy and, therefore, beyond~the scope of the present work. However, the~detailed derivation of  identity \eqref{TTMF1}, which is available in~\cite{Abrarov2017}, can be briefly outlined as a determination of the value $\eta$ at a given $k$ and $\gamma$ in the two-term Machin-like formula of kind
\begin{equation}\label{TTMF2}
\frac{\pi}{4} = 2^{k - 1}\arctan\left(\frac{1}{\gamma}\right) + \arctan\left(\frac{1}{\eta}\right)
\end{equation}
and a subsequent reformulation of Equation~\eqref{TTMF2} into Equation~\eqref{TTMF1} with the help of de~Moivre's~formula
$$
\left(\cos(x) + i\sin(x)\right)^n = \cos(nx) + i\sin(nx).
$$

It is not difficult to prove that if the first constant $\gamma$ of Equation \eqref{TTMF2} is an integer or a rational number, then the second constant $\eta$ must also be a rational number. Specifically, when $\gamma$ is either an integer or a rational number, then from the relation~\cite{Abrarov2017}
\begin{equation}\label{EI} 
\eta = \frac{2}{\left(\frac{\gamma + i}{\gamma - i}\right)^{2^{k-1}} - i} -i
\end{equation}
it follows that both the~real and the imaginary parts of $\eta $ must be rational numbers. However, from~the identity
$$
\eta = \frac{\cos\left(2^{k-1}\arctan\left(\frac{2\gamma}{\gamma^2 - 1}\right)\right)}{1 - \sin\left(2^{k-1}\arctan\left(\frac{2\gamma}{\gamma^2 - 1}\right)\right)}
$$
it follows that $\eta\in\mathbb{R}$ since $\gamma\in\mathbb{R}$. Consequently, we prove that $\eta$ is a rational number at $\gamma\in\mathbb{Q}$.

It should be noted that the two-term Machin-like formula of kind \eqref{TTMF2} that we considered in our paper~\cite{Abrarov2017} represents a practical interest. Recently, a~group of independent researchers has developed a different method for determination of the value of $\eta$ at given values of $k$ and $\gamma$ (see Table~2 in~\cite{Gasull2023}). Unlike our iterative method described in~\cite{Abrarov2017}, their algorithm is built on the basis of the rational functions $R_j(n,x)$ with a remarkable property (see~\cite{Gasull2023} a for detailed description).
$$
R_j(n,x) = \tan(n\theta + nj/4), \qquad x = \tan\theta, \qquad j\in\mathbb{N}.
$$

Following the same procedure that is described in our work~\cite{Abrarov2017}, one can also obtain a generalization of  identity \eqref{TTMF1} in the form
\begin{equation}\label{TTMF3}
\frac{\pi}{4} = \varphi\arctan\left(\frac{1}{\gamma}\right) + \arctan\left(\frac{1 - \sin\left(\varphi\arctan \left(\frac{2\gamma}{\gamma^2 - 1}\right)\right)}{\cos\left(\varphi\arctan\left(\frac{2\gamma}{\gamma^2 - 1}\right)\right)}\right).
\end{equation}
Consequently, the~two-term Machin-like Formula \eqref{TTMF2} can be generalized as
\begin{equation}\label{TTMF4}
\frac{\pi}{4} = \varphi\arctan\left(\frac{1}{\gamma}\right) + \arctan\left(\frac{1}{\eta}\right),
\end{equation}
where (compare with Equation \eqref{EI} above)
$$
\eta = \frac{2}{\left(\frac{\gamma + i}{\gamma - i}\right)^\varphi - i} - i.
$$

It should be noted that using the identity \cite{Abrarov2022}
\[
\arctan(x + y) = \arctan(x) + \arctan\left(\frac{y}{1 + (x + y)x}\right)
\]
the equation \eqref{TTMF4} can be represented as
\[
\frac{\pi}{4} = \arctan\left(\frac{\eta - 1}{\eta + 1}\right) + \arctan\left(\frac{1}{\eta}\right),
\]
from which it follows that
\[
\varphi\arctan\left(\frac{1}{\gamma}\right) = \arctan\left(\frac{\eta - 1}{\eta + 1}\right),
\]
where
\[
\frac{\eta - 1}{\eta + 1} = \frac{2i}{\left(\frac{\gamma + 1}{\gamma - 1}\right)^\varphi + 1} - i.
\]

Both constants $\varphi$ and $\gamma$ in Equation \eqref{TTMF3} may be chosen conveniently. For~example, consider a ratio of $22/7$, representing a rough approximation of $\pi$. Therefore, we can write
$$
\frac{22}{28} \approx 22\arctan\left(\frac{1}{28}\right) \approx \frac{\pi}{4}.
$$
Thus, by~choosing $\varphi = 22$ and $\gamma = 28$, we can show that the ratio
\[
\frac{1 - \sin\left(22\arctan\left(\frac{2 \times 28}{28^2 - 1}\right)\right)}{\cos\left(22\arctan\left(\frac{2 \times 28}{28^2 - 1}\right)\right)} = \frac{1744507482180328366854565127}{98646395734210062276153190241239}
\]
since
$$
\eta =\frac{2}{\left(\frac{28 + i}{28 - i}\right)^{22} - i} - i = \frac{98646395734210062276153190241239}{1744507482180328366854565127}.
$$
Consequently, we obtain the identity that corresponds to the first row of Table~1 in~\cite{Gasull2023}
\[
\frac{\pi}{4} = 22\arctan\left(\frac{1}{28}\right) + \arctan\left(\frac{\overbrace{1744507482180328366854565127}^{28\,{\rm digits}}}{\underbrace{98646395734210062276153190241239}_{32\,{\rm digits}}}\right).
\]
We can see now that Equation \eqref{TTMF1} is just a specific case $\varphi = 2^{k-1}$ of more general form of the two-term Machin-like Formula \eqref{TTMF3} for $\pi$. Therefore, the~method described in our paper~\cite{Abrarov2017} can also be generalized to generate all two-term Machin-like formulas of kind~\eqref{TTMF4}, shown in Table~1 from the work~\cite{Gasull2023} (see first row showing $f_{32}^{28}$).

In our recent publication~\cite{Abrarov2021a}, using the two-term Machin-like Formula \eqref{TTMF1} for $\pi$, we found the following iterative formula
\begin{equation}\label{IFSC1} 
\theta_n = \frac{1}{\frac{1}{\theta_{n-1}} + \frac{1}{2^k}\left(1 - \tan\left(\frac{2^{k-1}}{\theta_{n-1}}\right)\right)}, \qquad k \ge 1,
\end{equation}
where $\theta_1 = 2^{-k}$, such that
\begin{equation}\label{L1} 
\frac{\pi}{4} = 2^{k - 1}\lim_{n\to\infty}\,\frac{1}{\theta_n}.
\end{equation}
This equation can be employed to compute digits of the constant $\pi$ with quadratic convergence (see Mathematica code provided in~\cite{Abrarov2021a}).

Motivated by recent publications~\cite{Gasull2023,WolframCloud,Campbell2023,Maritz2023} in connection to our works~\cite{Abrarov2022,Abrarov2017,Abrarov2021a,Abrarov2018,Abrarov2021b}, we developed a new algorithm based on a modified version of the iterative Formula~\eqref{IFSC1}.

Although Equation \eqref{IFSC1} provides squared convergence in computing digits of $\pi$, its direct application results in a slow convergence rate in the intermediate steps of the calculations of the tangent function. This occurs because the argument of the tangent function in this equation tends to the relatively large value of $\pi/4$ as $n$ increases. In~this work, we propose a new method that resolves this problem. In~particular, we show how Equation \eqref{IFSC1} can be rearranged and used in combination with arctangent terms of the Machin-like formula for $\pi$. Such an approach may be promising for efficient computation of $\pi$ with more than $17$ digits per increment of $n$ in Equation \eqref{IF4T} that will be discussed below. To~the best of our knowledge, algorithmic implementation based on the combination of iterative and Machin-like formulas for computing digits of $\pi$ has never been reported.

\section{Preliminaries}

\subsection{Machin-Like~Formulas}

In our recent work, we  derived the following identity~\cite{Abrarov2022}
\small
\begin{equation}\label{BF} 
\frac{\pi}{4} = 2^{k - 1}\arctan\left(\frac{1}{A_k}\right) + \left(\sum\limits_{m = 1}^M \arctan \left( \frac{1}{\left\lfloor B_{m,k}\right\rfloor}\right)\right) + \arctan\left(\frac{1}{B_{M + 1,k}} \right),
\end{equation}
\normalsize
where
\begin{equation}\label{AkC} 
A_k = \left\lfloor \frac{a_k}{\sqrt{2 - a_{k - 1}}}\right\rfloor
\end{equation}
such that $a_0 = 0$ and $a_k = \sqrt{2 + a_{k - 1}}$ are nested radicals of $2$ and
\begin{equation}\label{BmkC} 
B_{m,k} = \frac{1 + \left\lfloor B_{m - 1,k} \right\rfloor B_{m - 1,k}}{\left\lfloor B_{m - 1,k} \right\rfloor - B_{m - 1,k}}, \qquad m \ge 2
\end{equation}
with an initial value $B_{1,k}$ that can be computed by substituting $A_k$ into Equation \eqref{EI}
\begin{equation}\label{B1kC} 
B_{1,k} = \frac{2}{\left(\frac{A_k + i}{A_k - i}\right)^{2^{k - 1}} - i} - i.
\end{equation}

Equation \eqref{BF} implies two important rules. First, since the integer $B_{0,k}$ is not defined, it follows that at $M = 0$, the sum of arctangent functions
\[
\left.\sum\limits_{m = 1}^M\arctan \left(\frac{1}{\left\lfloor B_{m,k} \right\rfloor}\right)\right|_{M = 0} = 0.
\]
Second, if~$\left\lfloor B_{M + 1,k} \right\rfloor - B_{M + 1,k} = 0$, then no further iteration is required, as~the fractional part of the number $B_{M + 1,k}$ does not~exist.

We may compute the coefficient $B_{1,k}$ by using Equation \eqref{B1kC} at smaller values of the integer $k$. However, the~computation slows down as $k$ increases. To~resolve this problem, we proposed a more efficient method of computation based on a two-step iterative formula~\cite{Abrarov2017}
\begin{equation}\label{TSIF} 
\left\{
\begin{aligned}
u_n &= u_{n - 1}^2 - v_{n - 1}^2, \\
v_n &= 2{u_{n - 1}}{v_{n - 1}}, \qquad n = 2,3,4,\ldots,k
\end{aligned}
\right.
\end{equation}
with initial values
\[
u_1 = \frac{A_k^2 - 1}{A_k^2 + 1}
\]
and
\[
v_1 = \frac{2A_k}{A_k^2 + 1}
\]
leading to
\begin{equation}\label{B1kF} 
B_{1,k} = \frac{u_k}{1 - v_k}.
\end{equation}

Equation \eqref{B1kF}, based on the two-step iteration \eqref{TSIF}, is more efficient for computation of the constant $B_{1,k}$ than Equation \eqref{B1kC}, since at larger values of the integer $k$, the rapidly growing exponent $2^{k - 1}$ in Equation \eqref{B1kC} drastically decelerates the~computation.

Consider a few examples. At~$k = 2$, we can find that
$$
A_2 = \left\lfloor\frac{\sqrt{2 + \sqrt{2}}}{\sqrt{2 - \sqrt{2}}}\right\rfloor = 2
$$
and at $M = 0$, the~constant $B_{1,2} = -7$ according to Equation \eqref{B1kF}. Since \scalebox{.95}[1.0]{$\left\lfloor B_{1,2} \right\rfloor - B_{1,2} = 0$}, the~constants $B_{m,2}$ at $m \ge 2$ do not exist. Consequently, we get the identity
\begin{equation}\label{HF} 
\begin{aligned}
\frac{\pi}{4} &= 2^{2 - 1}\arctan\left(\frac{1}{A_2}\right) + \arctan\left(\frac{1}{B_{1,2}}\right) \\ 
&= 2\arctan\left(\frac{1}{2}\right) - \arctan\left(\frac{1}{7}\right)
\end{aligned}
\end{equation}
that is commonly known as Hermann's formula~\cite{Borwein1987}.

At $k = 3$, we obtain
$$
A_3 = \left\lfloor\frac{\sqrt{2 + \sqrt{2 + \sqrt{2}}}}{\sqrt{2 - \sqrt{2 + \sqrt{2}}}}\right\rfloor = 5
$$
and at $M = 0$, the~coefficient $B_{1,3} = -239$ according to Equation \eqref{B1kF}. Since $\left\lfloor B_{1,3} \right\rfloor - B_{1,3} = 0$, the~constants $B_{m,3}$ at $m \ge 2$ do not exist. Consequently, we end up with the following identity
$$
\frac{\pi}{4} = 2^{3 - 1}\arctan\left(\frac{1}{A_3}\right) + \arctan\left(\frac{1}{B_{1,3}}\right)
$$
representing the original Machin-like Formula \eqref{OMF} for $\pi$.

The case $k = 4$ requires more computations. In~particular, we can see that
$$
A_4 = \left\lfloor\frac{\sqrt{2 + \sqrt{2 + \sqrt{2 + \sqrt{2}}}}}{\sqrt{2 - \sqrt{2 + \sqrt{2 + \sqrt{2}}}}}\right\rfloor = 10
$$
and at $M = 0$, we have
$$
B_{1,4} = -\frac{147153121}{1758719}.
$$
Consequently, we can write
$$
\begin{aligned}
\frac{\pi}{4} &= 2^{4-1}\arctan\left(\frac{1}{A_4}\right) + \arctan\left(\frac{1}{B_{1,4}}\right) \\ 
&= 8\arctan\left(\frac{1}{10}\right) - \arctan\left(\frac{1758719}{147153121}\right).
\end{aligned}
$$
Since $\left\lfloor B_{1,4} \right\rfloor - B_{1,4}\ne 0$, the~constant $B_{2,4}$ exists and can be computed according to Equation~\eqref{BmkC}.

Thus, at~$M = 1$, this leads to
$$
\frac{\pi}{4} = 8\arctan\left(\frac{1}{10}\right) - \arctan\left(\frac{1}{84}\right) - \arctan\left( \frac{579275}{12362620883}\right),
$$
where
$$
\frac{1}{84} = -\frac{1}{\left\lfloor B_{1,4}\right\rfloor}
$$
and
$$
\frac{579275}{12362620883}=-\frac{1}{{{B}_{2,4}}}.
$$

Again, since $\left\lfloor B_{2,4} \right\rfloor - B_{2,4} \ne 0$, the~constant  $B_{3,4}$ can be computed, and at $M = 2$, we can derive the following identity
$$
\begin{aligned}
\frac{\pi}{4} &= 8\arctan\left(\frac{1}{10}\right) - \arctan\left(\frac{1}{84}\right) - \arctan\left(\frac{1}{21342}\right) \\
&-\arctan\left(\frac{266167}{263843055464261}\right),
\end{aligned}
$$
where
$$
\frac{1}{21342} = -\frac{1}{\left\lfloor B_{2,4}\right\rfloor}
$$
and
$$
\frac{266167}{263843055464261} = -\frac{1}{B_{3,4}}.
$$

Repeating the same procedure over and over again up to $M = 5$, we can finally obtain the following seven-term Machin-like formula
\begin{equation}\label{STMF} 
\begin{aligned}
\frac{\pi}{4} &= 8\arctan\left(\frac{1}{10}\right) - \arctan\left(\frac{1}{84}\right) - \arctan\left(\frac{1}{21342}\right) \\
& -\arctan\left(\frac{1}{991268848}\right) - \arctan\left(\frac{1}{193018008592515208050}\right) \\
& -\arctan\left(\frac{1}{197967899896401851763240424238758988350338}\right) \\ 
& -\arctan\left(\frac{1}{\left| B_{6,4} \right|}\right),\\
\end{aligned}
\end{equation}
where
\[
\begin{aligned}
B_{6,4} &= -117573868168175352930277752844194126767991915008537\cdots \\
& 018836932014293678271636885792397.  
\end{aligned}
\]
is also an integer. Therefore, all iterations are~completed.

From these examples, we can see that Hermann's \eqref{HF}, Machin's \eqref{OMF} and the derived \eqref{STMF} formulas for $\pi$ belong to the same generic group, since all of them can be constructed from their generalized form \eqref{BF} at different integers $k$, equal to $2$, $3$ and $4$, respectively. Further, we will use Equation \eqref{STMF} as an example for computing digits of $\pi$.

The following Mathematica code:

\small
\begin{verbatim}
(* Define long string *)
longStr = StringJoin["11757386816817535293027775284419412676",
"7991915008537018836932014293678271636885792397"];

coeff = (10 + I)^8*(84 + I)^-1*(21342 + I)^-1*
(991268848 + I)^-1*(193018008592515208050 + I)^-1*
(197967899896401851763240424238758988350338 + I)^-1*
(FromDigits[longStr] + I)^-1;

Re[coeff] == Im[coeff]
\end{verbatim}
\normalsize
validates Equation \eqref{STMF} by returning {\ttfamily\bf True}. This code applies the product relation \eqref{PR} for~verification.

\subsection{Tangent~Function}

Since the tangent function can be represented as a series expansion
\begin{equation}\label{METF} 
\begin{aligned}
\tan(x) &= \sum\limits_{n = 1}^\infty \frac{(-1)^{n-1} 2^{2n}(2^{2n} - 1)B_{2n}}{(2n)!} x^{2n-1} \\ 
 &= x+\frac{x^{3}}{3}+\frac{2x^{5}}{15} + \frac{17{x^7}}{315} + \frac{62x^9}{2835} + \cdots,  
\end{aligned}
\end{equation}
where $B_{2n}$ are the Bernoulli numbers, defined by the contour integral
$$
B_n = \frac{n!}{2\pi i}\oint\frac{z}{e^z - 1}\frac{dz}{z^{n + 1}},
$$
it follows that
\begin{equation} \label{TBON} 
\tan(x) = x + O(x^3).
\end{equation}
Consequently, the accuracy of the tangent function improves with decreasing the argument $x$. 

Unfortunately, the~series expansion of the tangent function \eqref{METF} requires the determination of the Bernoulli numbers, which is itself a big challenge~\cite{Knuth1967,Harvey2010,Bailey2013,Beebe2017}. For~example, one of the most known formulas for computation of the Bernoulli numbers
\[
B_n = \sum\limits_{m = 0}^n\frac{1}{m + 1}\sum\limits_{\ell = 0}^m (-1)^\ell{m\choose\ell}\ell^n,
\]
is based on double summation with the binomial coefficients that decelerate the~computation.

There are other equations for computation of the tangent function~\cite{Bailey2013,Beebe2017,Havil2012} and one of the efficient techniques to perform computation of the tangent function is the Newton--Raphson iteration (see~\cite{Abrarov2021b} for more details). In~particular, the~following iteration formula
\begin{equation}\label{NRI} 
s_n(x) = s_{n - 1}(x) - \left(1 + s_{n - 1}^{2}(x) \right)\left(\arctan\left(s_{n - 1}(x)\right) - x\right)
\end{equation}
with an initial value that can be taken as
$$
s_{1}(x) = x,
$$
can be employed. This iteration leads to
$$
\tan(x) = \lim_{n\to\infty}\,s_n.
$$

Iteration \eqref{NRI} leads to the quadratic convergence of the tangent function. In~practice, however,~quadratic convergence can be achieved only if a sufficiently large number of the summation terms are applied in the series expansions like \eqref{MSE4AF}, \eqref{EF} or \eqref{AFSE1} to approximate the arctangent function. We have already applied the Newton--Raphson iteration \eqref{NRI}   to compute the digits of $\pi$ \cite{Abrarov2021b}.

As an option, we can also apply the most common equation
$$
\tan(x) = \frac{\sin(x)}{\cos(x)} = \frac{\sum\limits_{n = 0}^\infty (-1)^n\frac{x^{2n+1}}{(2n+1)!}}{\sum\limits_{n = 0}^\infty (-1)^n\frac{x^{2n}}{(2n)!}},
$$
where sine and cosine functions in the numerator and denominator are represented as the Maclaurin expansions. Unfortunately, this representation is not optimal for practical application since it requires separate computations of the expansion terms for the sine and cosine functions. However, if~we 
rewrite the tangent function as
$$
\tan(x) = \frac{\sin(x)}{\cos(x)}\,\frac{\sin(2x)}{\sin(2x)}
$$
and take into account that $\sin(2x) = 2\sin(x)\cos(x)$, we can obtain the following identity
\begin{equation}\label{TFR} 
\tan(x) = \frac{2\sin^2(x)}{\sin(2x)}.
\end{equation}

Although this this representation of the tangent function is not common, its application is significantly advantageous, since each $n$th term in the expansion
$$
\sin(x) = x - \frac{x^3}{3!} + \frac{x^5}{5!} - \frac{x^7}{7!} + \cdots\,,
$$
can be utilized again just by multiplying ${{2}^{2n+1}}$ to obtain the expansion for the denominator
\[
\sin(2x) = (2)x - (2^3)\frac{x^3}{3!} + (2^5)\frac{x^5}{5!} - (2^7)\frac{x^7}{7!} + \cdots \,.
\]
Such a technique may accelerate computation and reduce memory usage. Thus, we can write the following iterative procedure
\[
{{p}_{0}}\left( x \right)=0,
\]
\[
{{q}_{0}}\left( x \right)=0,
\]
\[
{{r}_{n}}\left( x \right)=\frac{{{\left( -1 \right)}^{n}}{{x}^{2n+1}}}{\left( 2n+1 \right)!},
\]
\[
{{p}_{n}}\left( x \right)={{p}_{n-1}}\left( x \right)+{{r}_{n-1}}\left( x \right),
\]
\[
{{q}_{n}}\left( x \right)={{q}_{n-1}}\left( x \right)+{{2}^{2n-1}}{{r}_{n-1}}\left( x \right),
\]
leading to
\begin{equation}\label{IF4T} 
\tan(x) = \lim_{n\to\infty}\,\frac{2p_n^2(x)}{q_n(x)},
\end{equation}
according to Equation \eqref{TFR}.

In this work, we used truncated Equation \eqref{IF4T} as an alternative to Equations \eqref{METF}~and~\eqref{NRI} since one of our objectives in this work is to develop an algorithm for computing the digits of $\pi$ as simply as possible and without undesirable surd~numbers.

\subsection{Arctangent~Function}

Since in this work we utilize the terms from the Machin-like Formula \eqref{BF}, a~series expansion with rapid convergence of the arctangent function should be applied. Based on our empirical results, we can consider two equations with rapid convergence that can be used for implementation. The~first equation is Euler's series expansion~\cite{Castellanos1988,Chien-Lih2005}
\begin{equation}\label{EF} 
\arctan(x) = \sum\limits_{n = 0}^\infty \frac{2^{2n}(n!)^2}{(2n+1)!}\frac{x^{2n + 1}}{\left(1 + x^2 \right)^{n+1}}.
\end{equation}

The second equation is given by the series expansion~\cite{Abrarov2023}
\small
\begin{equation}\label{AFSE1} 
\arctan(x) = 2\sum\limits_{m = 1}^M \sum\limits_{n = 1}^\infty \frac{1}{(2n - 1)(2m - 1)^{2n-1}}\,\frac{\kappa_n(x,\gamma_{m,M})}{\kappa_n^2(x,\gamma_{m,M}) + \lambda_n^2(x\gamma_{m,M})},
\end{equation}
\normalsize
where the expansion coefficients are computed by a two-step iteration such that
$$
\kappa_1(x,t) = 1/(xt),
$$
$$
\lambda_1(x,t) = 1,
$$
$$
\kappa_n(x,t) = \kappa_{n-1}(x,t)\left(1 - 1/(xt)^2\right) + 2\lambda_{n - 1}(x,t)/(xt),
$$
$$
\lambda_n(x,t) = \lambda_{n-1}(x,t)\left(1 - 1/(xt)^2\right) - 2\kappa_{n - 1}(x,t)/(xt),
$$
and
$$
\gamma_{m,M} = \frac{m - 1/2}{M}.
$$

The derivation of the series expansion \eqref{AFSE1}, shown in our work~\cite{Abrarov2023}, is based on the Enhanced Midpoint Integration (EMI) formula (see also~\cite{Abrarov2017})
\[
\int\limits_0^1 f(x,t)dt = 2\sum\limits_{m = 1}^M \sum\limits_{n = 0}^\infty \frac{1}{(2M)^{2n + 1}(2n + 1)!}\left.\frac{\partial^{2n}}{\partial t^{2n}}f(x,t)\right|_{t = \frac{m - 1/2}{M}},
\]
where we imply that
\[
f(x,t) = \frac{x}{2}\left(\frac{1}{1 + ixt} + \frac{1}{1 - ixt}\right)
\]
since
\[
\arctan(x) = \int\limits_0^1 \frac{x}{1 + x^2 t^2}dt = \int\limits_0^1 \frac{x}{2}\left(\frac{1}{1 + ixt} + \frac{1}{1 - ixt}\right)dt.
\]

The series expansion \eqref{AFSE1} appears to be rapid in convergence. In~particular, even at $M=1$, its convergence rate is faster by many orders of the magnitude than that of Euler's Formula \eqref{EF} (see~\cite{Abrarov2023} for details). Thus, by~taking $M = 1$, Equation \eqref{AFSE1} can be conveniently rearranged as
\begin{equation}\label{AFSE2} 
\arctan(x) = 2\sum\limits_{n = 1}^\infty \frac{1}{2n - 1}\frac{g_n(x)}{g_n^2(x) + h_n^2(x)},
\end{equation}
where the expansion coefficient can be computed by the following two-step iteration
\[
g_n(x) = g_{n-1}(x)\left(1 - 4/x^2\right) + 4h_{n - 1}(x)/x,
\]
\[
h_n(x) = h_{n - 1}(x)\left(1 - 4/x^2\right) - 4g_{n-1}(x)/x,
\]
with initial values 
\[
g_1(x) = 2/x,
\]
\[
h_1(x) = 1.
\]

Series expansions \eqref{EF} and \eqref{AFSE2} are significantly faster in convergence than the Maclaurin series expansion \eqref{MSE4AF}. In~this work, the~truncated series expansion \eqref{AFSE2} is used. Although a~value of  integer $M$ that is greater than $1$  further improves the convergence rate, it may be preferable to apply Equation \eqref{AFSE2} rather than its generalization \eqref{AFSE1}, since an increment of $M$ just by $1$ increases the number of the terms in  series expansion \eqref{AFSE1} by a factor of $2$.

\section{Results and~Discussion}

\subsection{Modified~Iteration}

Changing the variable $\theta_n \to 1/\sigma_n$ in Equation \eqref{IFSC1} leads to a more convenient form
\begin{equation}\label{IFSC2} 
\sigma_n = \sigma_{n - 1} + 2^{-k}\left(1-\tan\left(2^{k - 1}\sigma_{n - 1}\right)\right), \qquad k \ge 1,
\end{equation}
where
\[
\sigma_1 = 2^{-k}.
\]
Consequently, the~constant $\pi$ can be found by iteration in accordance with Equation \eqref{L1} such that
\begin{equation}\label{L2}
\frac{\pi}{4} = 2^{k - 1}\lim_{n\to\infty}\,\sigma_n.
\end{equation}

Comparing the limit \eqref{L2} with (see~\cite{Abrarov2018} for derivation)
$$
\frac{\pi}{4} = 2^{k - 1}\arctan\left(\frac{\sqrt{2 - a_{k - 1}}}{a_k}\right),
$$
we can obtain an important relation
$$
\sigma_n\to\arctan\left(\frac{\sqrt{2 - a_{k - 1}}}{a_k}\right), \qquad n \to \infty
$$
or
\begin{equation}\label{R1} 
2^{k - 1}{\sigma_n}\to\frac{\pi}{4}, \qquad n \to \infty.
\end{equation}
Consequently, the~argument of the tangent function in Equation \eqref{IFSC2} tends to $1$ with increasing~$n$.

Since, in Equation \eqref{IFSC2}, the argument of the tangent function tends to $1$, the~convergence per iteration of $n$ in applied Equation \eqref{IF4T} is expected to be extremely slow. Indeed, if~the argument of the tangent function is not small enough, then, in accordance with relation \eqref{TBON}, we cannot gain a reasonable convergence of the tangent~function.

This problem can be effectively resolved by introducing a constant such that
$$
c \approx \frac{1}{2^{k - 1}}\frac{\pi}{4}.
$$
The value of this constant is close to $\sigma_n$ when $n$ is sufficiently large (see Equation \eqref{R1}). Therefore, with~the help of this constant, we can modify the iteration Formula \eqref{IFSC2} by using the following procedure
\[
\begin{aligned}
\sigma_1 &= 2^{-k}, \\
\delta_1 &= c - \sigma_1,\\
\sigma_2 &= c - \delta_1 + 2^{-k}\left(1 -\tan\left(2^{k - 1}(c - \delta_1)\right)\right), \\
\delta_2 &= c - \sigma_2, \\
\sigma_3 &= c - \delta_2 + 2^{-k}\left(1 - \tan\left(2^{k - 1}(c - \delta_2)\right)\right), \\
\delta_3 &= c - \sigma_3, \\
&\vdots
\end{aligned}
\]
\[
\begin{aligned}
&\hspace{1cm}\sigma_n = c - \delta_{n - 1} + 2^{-k}\left(1 - \tan\left(2^{k - 1}(c - \delta_{n - 1})\right)\right), \\
&\hspace{1cm}\delta_n = c - \sigma_n.
\end{aligned}
\]

The expression for $\sigma_n$ can be simplified as follows
\[
\begin{aligned}
\sigma_n &= c - \delta_{n - 1} + 2^{-k}\left(1 - \tan\left(2^{k-1}(c - \delta_{n - 1})\right)\right) \\
&= c - \delta_{n - 1} + 2^{-k}\left(1 - \tan\left(2^{k-1} c - 2^{k - 1}\delta_{n - 1}\right)\right) \\ 
&= c - \delta_{n - 1} + 2^{-k}\left(1 - \frac{\tan\left(2^{k - 1}c\right) - \tan\left(2^{k - 1}\delta_{n - 1}\right)}{1 + \tan\left(2^{k - 1} c\right)\tan\left(2^{k - 1}\delta_{n - 1}\right)}\right)  
\end{aligned}
\]
or
\[
\sigma_n = c - \delta_{n - 1} +  2^{-k}\left(1 -\frac{\alpha - \tan\left(2^{k - 1}\delta_{n-1}\right)}{1 + \alpha\tan\left(2^{k - 1}\delta_{n - 1}\right)}\right)
\]
or
\begin{equation}\label{IFSC3} 
\sigma_n = \sigma_{n - 1} + 2^{-k}\left(1 - \frac{\alpha - \tan\left(2^{k - 1}\delta_{n - 1}\right)}{1 +\alpha\tan\left(2^{k - 1}\delta_{n - 1}\right)}\right),
\end{equation}
where $\alpha = \tan\left(2^{k - 1}c\right)$. Since $\left|\delta_n\right|$ is supposed to be a small value, we may expect a reasonable convergence rate in  iteration \eqref{IFSC3} at each increment $n$ in Equation \eqref{IF4T}.

The tangent function in Equation \eqref{IFSC3} is represented twice. Therefore, we can define
\[
\tau_n = \tan\left(2^{k - 1}\delta_n\right)
\]
and rewrite this equation in a more simplified form
\begin{equation}\label{IFSC4} 
\sigma_n = \sigma_{n - 1} + 2^{-k}\left(1 - \frac{\alpha - \tau_{n - 1}}{1 + \alpha\tau_{n-1}}\right).
\end{equation}

\subsection{Methodology}

In our algorithm, we utilize a modified Equation \eqref{IFSC4}. This yields the iteration based on a set of equations
\begin{equation}\label{SE} 
\left\{
\begin{aligned}
& \delta_{n - 1} = c - \sigma_{n - 1}, \\ 
& \tau_{n - 1} = \tan\left(2^{k - 1}\delta_{n - 1}\right), \\ 
& \sigma_n = \sigma_{n - 1} + 2^{-k}\left(1 - \frac{\alpha - \tau_{n - 1}}{1 + \alpha\tau_{n - 1}}\right).
\end{aligned}
\right.
\end{equation}
We can take a few initial terms in the Machin-like Formula \eqref{BF} as a value to compute the constant $c$.

Return to Equation \eqref{STMF}, corresponding to the case $k = 4$. If~we take only the first term, then the constant
\begin{equation}\label{CC1} 
c = \frac{1}{2^{k - 1}} \,\, \underbrace{8\arctan\left(\frac{1}{10}\right)}_{\rm first \, term \, of~\, Equation \, (24)} = \arctan\left(\frac{1}{10}\right)
\end{equation}
is not close enough to $\pi/(4 \times 2^{k - 1})$ to compute the digits of $\pi$. As~a consequence, we cannot achieve a rapid convergence. Specifically, our empirical results show that with constant \eqref{CC1}, the~set \eqref{SE} of iteration formulas provides convergence of four to five digits of $\pi$ per increment of $n$ in Equation \eqref{IF4T}. However, if~we take the first two terms from Equation \eqref{STMF} such that
\begin{equation}\label{CC2} 
\begin{aligned}
c &= \frac{1}{2^{k - 1}}\left(\underbrace{8\arctan\left(\frac{1}{10}\right)}_{\rm first \, term \, of~\, Equation \, (24)} \, \underbrace{-\arctan\left(\frac{1}{84}\right)}_{\rm second \, term \, of~\, Equation \, (24)}\right) \\
&= \arctan\left(\frac{1}{10}\right) - \frac{1}{8}\arctan\left(\frac{1}{84}\right),
\end{aligned}
\end{equation}
then the convergence rate significantly improves, providing ten digits of $\pi$ per~increment.

Consider the following identity~\cite{Oliver2012}
\begin{equation}\label{TI} 
\begin{aligned}
\tan(nx) &= \frac{1}{i}\,\,\frac{\left(1 + i\tan(x)\right)^n - \left(1 - i\tan(x) \right)^n}{\left(1 + i\tan(x)\right)^n + \left(1 - i\tan(x)\right)^n} \\
&= \frac{2i\left(1 - i\tan(x)\right)^n}{\left(1 - i\tan(x)\right)^n + \left(1 + i\tan(x)\right)^n} - i.
\end{aligned}
\end{equation}
This identity can be used for computation of the constant $\alpha$. Thus, according to identity~\eqref{TI}, we have
\[
\tan \left( 8\arctan \left( \frac{1}{10} \right) \right)=\frac{2i{{\left( 1-i/10 \right)}^{8}}}{{{\left( 1-i/10 \right)}^{8}}+{{\left( 1+i/10 \right)}^{8}}}-i=\frac{\text{74455920}}{\text{72697201}}.
\]
Now, using
\[
\tan\left(\arctan\left(\frac{1}{84}\right)\right) = \frac{1}{84}
\]
and an elementary trigonometric relation
\begin{equation}\label{ETR} 
\tan(x - y) = \frac{\tan(x) - \tan(y)}{1 + \tan(x)\tan(y)}
\end{equation}
we can find that the constant is
\begin{equation}\label{R4A} 
\alpha =\tan\left(8\arctan\left(\frac{1}{10}\right) - \arctan\left(\frac{1}{84}\right)\right) = \frac{6181600079}{6181020804}.
\end{equation}
	
Once the exact value of the constant $\alpha$ is calculated according to Equation \eqref{R4A}, we can use the set of iteration Formula \eqref{SE} for computing the digits of $\pi$. The~Mathematica codes and their description are provided in the next~section.

Although the convergence rate of $10$ digits of $\pi$ per increment of $n$ in Equation \eqref{IF4T} is relatively high, application of the Machin-like Formula \eqref{BF} may not be efficient at a smaller value of the integer $k$. In~particular, the~computation of both arctangents $\arctan(1/10)$ and $\arctan(1/84)$ is expected to be slow, since the two integers in the argument denominators ($10$ and $84$) are not big enough for rapid convergence in accordance with Equation \eqref{ABON}.

Consider another case where $k$ is sufficiently large. We can take, for~example, $k = 27$. Applying identity \eqref{BF} together with Equation \eqref{B1kF}, we can find that~\cite{Abrarov2017}
\small
\begin{equation}\label{TTMFk27} 
\begin{aligned}
\frac{\pi}{4} &= 2^{27 - 1}\arctan\left(\frac{1}{A_{27}}\right) + \arctan\left(\frac{1}{B_{1,27}}\right) 
\\ 
&= 67108864\arctan\left(\frac{1}{85445659}\right) - \arctan\left(\frac{\overbrace{9732933578\ldots 4975692799}^{522,185,807\,{\rm digits}}}{\underbrace{2368557598\ldots 9903554561}_{522,185,816\,\,{\rm digits}}}\right).
\end{aligned}
\end{equation}
\normalsize
Recently, the~same equation was derived by Gasull~et~al. by a different method of computation (see Table~2 in~\cite{Gasull2023} showing the row with $f_{522,185,816}^{522,185,807}$).

At first glance, it may appear problematic to find the corresponding coefficient $\alpha$, since at $k = 27$, the substitution $2^{27 - 1} = 67108864$ into identity \eqref{TI} results in an expression
\[
\begin{aligned}
\alpha &= \tan\left(67108864\arctan\left(\frac{1}{85445659}\right)\right) \\ 
&= \frac{2i\left(1 - i/85445659\right)^{67108864}}{\left(1 - i/85445659\right)^{67108864} + \left(1 + i/85445659\right)^{67108864}}-i  
\end{aligned}
\]
that is impossible to compute due to the extremely large value of the exponent. However, application of the two-term Machin-like formula for $\pi$ of kind \eqref{TTMF2} gives a big advantage, since the multiplier $2^{k - 1}$ is continuously divisible by $2$. Thus, using the elementary trigonometric identity
$$
\tan(2x) = \frac{2\tan(x)}{1 - \tan^2(x)}
$$
in the iteration process, we can compute the required constant $\alpha$. Specifically, in~accordance with this identity and due to relation $x = \tan\left(\arctan(x )\right)$, at~$k = 27$, the following iteration formula
\begin{equation}\label{lambda}
\lambda_n(x) = \frac{2\lambda_{n - 1}(x)}{1 - \lambda_{n - 1}^2(x)},
\end{equation}
with initial value
\[
\lambda_1(x) = \frac{2x}{1 - x^2},
\]
yields a rational number
\begin{equation}\label{alpha}
\begin{aligned}
\alpha &= \lambda_{27 - 1}\left(\frac{1}{85445659}\right) = \frac{\overbrace{1184278804\ldots 8037709539}^{522,185,816\,{\rm digits}}}{\underbrace{1184278794\ldots 2617464027}_{522,185,816\,{\rm digits}}} \\ 
&= 1.00000000821844790606\ldots,
\end{aligned}
\end{equation}
in which the integer $85445659 = \left\lfloor A_{27} \right\rfloor$ is calculated according to Equation \eqref{AkC}. It is interesting to note that the number of digits (522,185,816) in this equation is the same as the number of digits in the denominator of Equation \eqref{TTMFk27}.

At a larger value of the integer $k$, it is sufficient to take only the first term (leading term) of the Machin-like formula for $\pi$ to achieve  rapid convergence. In~particular, using only the leading term from Equation \eqref{TTMFk27}
$$
67108864\arctan\left(\frac{1}{85445659}\right),
$$
the set of iteration Formula \eqref{SE} provides $17$ to $18$ digits of $\pi$ per increment of $n$ in Equation~\eqref{IF4T}. The~corresponding Mathematica code and its description are discussed in the next section.

\section{Mathematica Codes and~Description}

The Mathematica codes consist of seven cells that can be copied and pasted directly to the Mathematica notebook. The~first cell is given by the following code:
\small
\begin{verbatim}
atanF := {(* ARCTANGENT FUNCTION APPROXIMATION *)

Clear[atan,g,h];

(* Expansion coefficients *)
g[1,x_] := g[1,x] = 2/x; 
h[1,x_] := h[1,x] = 1; 
g[m_,x_] := g[m,x] = g[m - 1,x]*(1 - 4/x^2) + 4*(h[m - 1,x]/x); 
h[m_,x_] := h[m,x] = h[m - 1,x]*(1 - 4/x^2) - 4*(g[m - 1,x]/x); 

(* Arctangent approximation *)
atan[x_,n_] := atan[x,n] = 2*Sum[(1/(2*m - 1))*(g[m,x]/
  (g[m,x]^2 + h[m,x]^2)),{m,1,n}]};
\end{verbatim}
\normalsize
that defines the two-step iterative method for the arctangent function, based on Equation~\eqref{AFSE2}. This function can be invoked by running the command {\ttfamily\bf atanF}.

The second cell provides the code:
\newpage
\small
\begin{verbatim}
tanF:= {(*TANGENT FUNCTION APPROXIMATION*)

Clear[p,q,r,t];

(* Computing coefficients *)
r[n_,x_] := r[n,x] = (-1)^n*(x^(2*n + 1)/(2*n + 1)!);
p[0,x_] := p[0,x] = 0;
q[0, x_] := q[0,x] = 0;
p[n_,x_] := p[n,x] = p[n - 1,x] + r[n - 1,x];
q[n_,x_] := q[n,x] = q[n - 1,x] + 2^(2*n - 1)*r[n - 1,x];

(* Tangent approximation *)
t[x_,n_] := t[x,n] = 2*(p[n,x]^2/q[n,x])};
\end{verbatim}
\normalsize
that defines the iterative method for the tangent function, based on Equation \eqref{IF4T}. This function can be invoked by running the command {\ttfamily\bf tanF}.

The code below in the third cell:
\small
\begin{verbatim}
heading := {Print[Abs[MantissaExponent[Pi - 2^(k + 1)*
  \[Sigma][1]][[2]]]," digits of \[Pi] before iteration"];
    Print["-------------------------------"];
      Print["Number of terms n"," | ","Digits of \[Pi]"];
        Print["-------------------------------"]};

ending:={Print["-------------------------------"];
  Print[Abs[MantissaExponent[Pi - 2^(k + 1)*\[Sigma][2]][[2]]],
    " digits of \[Pi] after iteration"]};
\end{verbatim}
\normalsize
defines the output format that includes the heading and ending parts for intermediate computed data. The~header and the ending parts can be invoked by running the commands heading and ending, respectively.

\pagestyle{empty}
The code in the fourth cell:
\small
\begin{verbatim}
Clear[k,f,\[Alpha],c,\[Sigma],\[Delta],\[Tau],n];
(* Computing coefficient alpha *)
k=4; f[x_,n_] = (2*I*(1 - I*x)^n)/((1 - I*x)^n + (1 + I*x)^n) - I;
\[Alpha] = f[1/10,2^(k - 1)];
(* Computing coefficient c *)
atanF; c = atan[1/10,500]; 

(* Iteration *)
\[Sigma][1] = SetPrecision[Floor[Pi*10^100]/10^100,250]/2^(k + 1); 
\[Delta] = c - \[Sigma][1]; 
n = 1;
heading;
While[n <= 42, 
  tanF; \[Tau] = t[2^(k - 1)*\[Delta],n]; \[Sigma][2] =
    SetPrecision[\[Sigma][1] + (1 - (\[Alpha] - \[Tau])/
      (1 + \[Alpha]*\[Tau]))/2^k,5 + 5*n];
        piApp1 = 2^(k + 1)*\[Sigma][2]; 
          str=If[n < 10,"                 | ","                | "]; 
            If[n <= 5 || n >= 33, 
              Print[n,str,Abs[MantissaExponent[Pi - piApp1][[2]]]]]; 
                If[n==6,Print["...   COMPUTING     ..."]]; n++]; 
ending;
\end{verbatim}
\normalsize
computes the coefficient $\alpha$ according to Equation \eqref{TI} and coefficient $c$ at $k = 4$ by taking only the leading (first) term of the Machin-like Formula \eqref{STMF} for $\pi$. As~we do not require the highest precision at each cycle of computation within the while loop, the~precision increases with increasing $n$ in Equation \eqref{IF4T}. The~parameter of the precision is given as $5 + 5n$, where $n$ is the increment in Equation \eqref{IF4T}, since at $n=1$, the number of correct digits of $\pi$ is $5$ and multiplier $5$ is the largest number of digits of $\pi$ per increment $n$.

\pagestyle{plain}
Suppose we know the first $100$ digits of $\pi$. These digits of $\pi$ can be extracted from Mathematica by using the following division
\[
\frac{\left\lfloor 10^{100} \pi\right\rfloor}{10^{100}}.
\]
Once we know $100$ digits of $\pi$, we can use the described iteration to double its~number.

Running this cell generates the following output:
\small
\begin{verbatim}
100 digits of π before iteration
-------------------------------
Number of terms n | Digits of π
-------------------------------
1                 | 5
2                 | 9
3                 | 14
4                 | 19
5                 | 25
...   COMPUTING     ...
33                | 169
34                | 174
35                | 179
36                | 184
37                | 189
38                | 194
39                | 199
40                | 200
41                | 200
42                | 200
-------------------------------
200 digits of π after iteration
\end{verbatim}
\normalsize
that shows the first five and last fifteen cycles. As~we can see, each increment of $n$ by one gives four to five digits of $\pi$.	After~$39$ cycles, the~convergence slows down due to saturation. In~particular, after~40th cycle, the~number of digits  doubles to $200$. This saturation occurs since we have reached the limit of squared convergence in the determination of the digits of $\pi$.

\pagestyle{empty}
The code in the fifth cell:
\small
\begin{verbatim}
Clear[k,f,\[Alpha],c,\[Sigma],\[Delta],\[Tau],n];
(* Computing coefficient alpha *)
k=4; f[x_,n_] = (2*I*(1 - I*x)^n)/((1 - I*x)^n + (1 + I*x)^n) - I;
(* Computing coefficient c *)
\[Alpha] = f[1/10,2^(k - 1)]; \[Alpha] = 
  (\[Alpha] - 1/84)/(1 + \[Alpha]/84); 

atanF; c = atan[1/10,500] - (1/2^(k - 1))*atan[1/84,500];

(* Iteration *)
\[Sigma][1] = SetPrecision[piApp1/2^(k + 1),500];
\[Delta] = c-\[Sigma][1];
n = 1;
heading;
While[n <= 42,
  tanF; \[Tau] = t[2^(k - 1)*\[Delta],n]; \[Sigma][2] =
    SetPrecision[\[Sigma][1] + (1 - (\[Alpha] - \[Tau])/
      (1 + \[Alpha]*\[Tau]))/2^k,12 + 10*n];
        piApp2 = 2^(k + 1)*\[Sigma][2];
          str=If[n < 10,"                 | ","                | "];
            If[n <=5 || n >= 33,
              Print[n,str,Abs[MantissaExponent[Pi - piApp2][[2]]]]];
                If[n == 6,Print["...   COMPUTING     ..."]]; n++];
ending;
\end{verbatim}
\normalsize
computes the coefficients $c$ and $\alpha$ at $k = 4$ with the help of Equations \eqref{CC2}, \eqref{lambda} and \eqref{alpha}. Again, the~precision increases with increasing $n$ in Equation \eqref{IF4T}. The~parameter of the precision is given as $12 + 10n$, where $n$ is the increment in Equation \eqref{IF4T}, since at $n = 1$, the number of correct digits of $\pi$ is $12$ and multiplier $10$ is the largest number of digits of $\pi$ per increment $n$. We can use the $200$ already obtained digits of $\pi$ to double it by~iteration.

\pagestyle{plain}
Running this cell produces the output:
\small
\begin{verbatim}
200 digits of π before iteration
-------------------------------
Number of terms n | Digits of π
-------------------------------
1                 | 12
2                 | 21
3                 | 31
4                 | 41
5                 | 51
...   COMPUTING     ...
33                | 341
34                | 351
35                | 361
36                | 371
37                | 381
38                | 391
39                | 401
40                | 402
41                | 402
42                | 402
-------------------------------
402 digits of π after iteration
\end{verbatim}
\normalsize
that shows first five and last fifteen cycles. The~convergence rate always remains the same, $10$ digits per increment of $n$. After~$39$ cycles, the~convergence slows down. In~particular, after~the 40th cycle, the~number of digits of $\pi$  doubles to $402$. Further cycles do not contribute to increasing the number of digits since again we reached the limit in the determination of the digits of $\pi$.

The code in the sixth cell is the most interesting:
\small
\begin{verbatim}
Clear[k,\[Alpha],c,\[Sigma],\[Delta],\[Tau],n];
(* Computing coefficient alpha *)
k = 27; \[Sigma][1] = SetPrecision[piApp2/2^(k + 1),1000];
(* Computing coefficient c *)
c = SetPrecision[atan[1/85445659,500],1000];

(* Iteration *)
\[Alpha] = t[2^(k - 1)*c,500];
\[Delta] = c - \[Sigma][1];
n = 1;
heading;
While[n <= 46,
  tanF; \[Tau] = t[2^(k - 1)*\[Delta],n]; \[Sigma][2] =
    SetPrecision[\[Sigma][1] + (1 - (\[Alpha] - \[Tau])/
      (1 + \[Alpha]*\[Tau]))/2^k,25 + 18*n]; piApp3 = 2^(k + 1)*
			  \[Sigma][2];
          str=If[n < 10,"                 | ","                | "];
            If[n <=5 || n >= 37,
              Print[n,str,Abs[MantissaExponent[Pi - piApp3][[2]]]]];
                If[n == 6,Print["...   COMPUTING     ..."]]; n++];
ending;
\end{verbatim}
\normalsize
as it provides an excellent convergence rate at $17$ to $18$ digits per increment $n$ in Equation \eqref{IF4T} at $k = 27$. Although~the exact value of the rational number $\alpha$ can be used (see Equation~\eqref{alpha}), this code does not utilize it since a regular desktop or laptop computer requires a few hours for computation by using Equation \eqref{alpha}. However, we can observe a rapid convergence  computing this coefficient with tangent and arctangent approximations based on Equations~\eqref{alpha}~and~\eqref{AFSE2}, respectively. The~code in this cell utilizes only a leading term of the Machin-like Formula \eqref{TTMFk27} for $\pi$. This is possible to achieve since at larger value of $k = 27$, the~argument $1/85445659$ in the leading arctangent term becomes small. Consequently, this reduces the value of $\delta_n$ and improves the convergence rate according to Equation \eqref{IF4T}.

The precision of this code increases with increasing $n$ in Equation \eqref{IF4T}. The~parameter of the precision is given as $25 + 18n$, where $n$ is the increment in Equation \eqref{IF4T}, since at $n = 1$, the number of correct digits of $\pi$ is $25$ and the multiplier $18$ is the largest number of digits of $\pi$ per increment $n$. Again, we can use the $402$ already obtained  digits of $\pi$ to double it by~iteration.

Running this cell produces the output:
\small
\begin{verbatim}
402 digits of π before iteration
-------------------------------
Number of terms n | Digits of π
-------------------------------
1                 | 25
2                 | 42
3                 | 60
4                 | 78
5                 | 96
...   COMPUTING     ...
37                | 690
38                | 708
39                | 726
40                | 744
41                | 762
42                | 780
43                | 798
44                | 804
45                | 804
46                | 804
-------------------------------
804 digits of π after iteration
\end{verbatim}
\normalsize
that shows first five and last fifteen cycles. After~$43$ cycles, the~convergence slows down. In~particular, after the~44th cycle, the~number of digits of $\pi$  doubles to $804$. Further cycles do not contribute to the number of digits since we have reached the limit in the determination of the digits of $\pi$ for squared~convergence.

Since the proposed method requires at least the leading term of the Machin-like Formula \eqref{TTMFk27} for $\pi$, we have to verify its convergence rate to ensure efficient computation. Thus, at~$k = 27$, the corresponding value of the leading term in Equation \eqref{TTMFk27} is
\begin{equation}\label{FTAFk27} 
\arctan\left(\frac{1}{\left\lfloor A_{27} \right\rfloor}\right) = \arctan\left(\frac{1}{85445659}\right).
\end{equation}

An additional cell below is the Mathematica code:
\small
\begin{verbatim}
Clear[arg,n]

(* Argument of arctangent *)
arg=SetPrecision[1/85445659,500];

Print["Increment of n"," | ","Correct digits"];
Print["---------------------------------"];
n = 1;
(* Convergence of arctangent *)
While[n <= 15,
  str=If[n < 10,"              | ","             | "];
    Print[n,str,
      Abs[MantissaExponent[ArcTan[arg] - atan[arg,n]][[2]]]]; n++];
\end{verbatim}
\normalsize
showing the convergence rate of the arctangent function. This code generates following~output:
\small
\begin{verbatim}
Increment of n | Correct digits
---------------------------------
1              | 24
2              | 41
3              | 58
4              | 74
5              | 91
6              | 107
7              | 124
8              | 140
9              | 157
10             | 173
11             | 190
12             | 206
13             | 223
14             | 239
15             | 256
\end{verbatim}
\normalsize

As we can see, the~code provides $16$ to $17$ correct digits of the arctangent function value~\eqref{FTAFk27} per increment of $n$ in Equation \eqref{AFSE2}. Furthermore, computing just a single arctangent term greatly minimizes Lehmer's measure \eqref{LM}. In~particular, the value of the measure of the arctangent term \eqref{FTAFk27} is only
\[
\mu = \frac{1}{\log_{10}{(85445659)}} \approx 0.126077.
\]

The increase of the convergence rate with increasing $k$ is due to a decrease of the argument of the tangent function. These results show that the set of equations \eqref{SE} can provide  rapid convergence when computing the digits of $\pi$.

\section{Rational~Numbers}

There is an alternative method for computing the digits of $\pi$. This method is based only on rational numbers and implemented without the iterative Formula \eqref{IF4T}.

As an example, consider the following approximation
\small
\begin{equation}\label{PA} 
\frac{\pi}{4} \approx 8\arctan\left(\frac{1}{10}\right) - \arctan\left(\frac{1}{84}\right) - \arctan\left(\frac{1}{21342}\right) - \arctan\left(\frac{1}{991268848}\right)
\end{equation}
\normalsize
that represents the first four terms in Equation \eqref{STMF}. This approximation can provide at most $19$ decimal digits of $\pi$. Consequently, we can define
\small
\[
\begin{aligned}
\sigma_1 &= \arctan\left(\frac{1}{10}\right) - \frac{1}{8}\left[\arctan\left(\frac{1}{84}\right) + \arctan\left(\frac{1}{21342}\right) + \arctan\left(\frac{1}{991268848}\right)\right] \\
&= 0.098174770424681038702605213693...
\end{aligned}
\]
\normalsize
such that
\[
\frac{\pi}{4} \approx 8\sigma_1
\]
according to relation \eqref{PA}.

Using the iterative Formula \eqref{lambda} together with identity \eqref{ETR}, we can find the exact value of the tangent as a rational number
\begin{equation}\label{VT} 
\tan(8\sigma_1) = \frac{26153940164285810690885}{26153940164285810690614}
\end{equation}
since all arguments of the arctangent terms in approximation \eqref{PA} are also rational numbers (integer reciprocals). Once the exact value of tangent \eqref{VT} is known,  substituting it into iterative Formula \eqref{IFSC2} results in
\[
\sigma_2 = \sigma_1 + \frac{1}{16}\left(1 - \frac{26153940164285810690885}{26153940164285810690614}\right) 
\]
such that the value $32\sigma_2$ doubles the number of the decimal digits of $\pi$ up to $39$.

The corresponding Mathematica code:
\small
\begin{verbatim}
(* Iterative formula (45) *)
\[Lambda][x_,1] := \[Lambda][x,1] = (2*x)/(1 - x^2);
\[Lambda][x_,n_] := \[Lambda][x,n] = (2\[Lambda][x,n - 1])/
  (1 - \[Lambda][x,n - 1]^2);

(* Integer k *)
k = 4;

(* Computation of rational number *)
rNum = \[Lambda][1/10,k - 1];
rNum = (rNum - 1/84)/(1 + rNum*1/84);
rNum = (rNum - 1/21342)/(1 + rNum*1/21342);
rNum = (rNum - 1/991268848)/(1 + rNum*1/991268848);

(* Before iteration *)
\[Sigma][1] = N[ArcTan[1/10] - 1/2^(k - 1)*
  (ArcTan[1/84] + ArcTan[1/21342] + ArcTan[1/991268848]),100];
(* After iteration *)
\[Sigma][2]=\[Sigma][1] + 1/2^k*(1 - rNum);

Print[MantissaExponent[Pi - 2^(k + 1)*\[Sigma][1]][[2]]//Abs,
  " digits of \[Pi] before iteration"];
Print[MantissaExponent[Pi - 2^(k + 1)*\[Sigma][2]][[2]]//Abs,
  " digits of \[Pi] after iteration"];
\end{verbatim}
\normalsize
produces the following output:
\small
\begin{verbatim}
19 digits of π before iteration
39 digits of π after iteration
\end{verbatim}
\normalsize
showing number of digits of $\pi$ before and after a single iteration. This example demonstrates that the combination of iteration Formula \eqref{IFSC2} with the Machin-like formula for $\pi$ of type \eqref{BF} can also be implemented for computing the digits of $\pi$ without any approximation formula of the tangent~function.

Since Equation \eqref{STMF} consists of only seven terms, we should not expect a high accuracy. However, increasing $k$ in Equation \eqref{BF} leads to a rapid increase in the number of  arctangent terms. Therefore, at~sufficiently large values of $k$, we may also achieve the required accuracy for computing the digits of $\pi$ by using this alternative method based on rational numbers.

\section{Conclusions}

A new iterative method for computing the digits of $\pi$ by argument reduction of the tangent function is developed. This method combines a modified iterative formula for $\pi$ with squared convergence and a leading arctangent term from the Machin-like Formula~\eqref{BF}. The~computational test shows that algorithmic implementation can provide more than $17$ digits of $\pi$ per increment $n$ in Equation~\eqref{IF4T}. This method requires no surd numbers, and with an arbitrarily large $k$, there is no upper limit for the convergence rate.

\section*{Acknowledgment}

This work was supported by National Research Council Canada, Thoth Technology Inc., York University and Epic College of Technology. The authors wish to thank the reviewers for their constructive comments and~recommendations.

\end{document}